\documentstyle{amsppt}
\magnification1200
\pagewidth{30pc}
\pageheight{47pc}
\binoppenalty=10000\relpenalty=10000
\TagsOnRight
\loadbold
\UseAMSsymbols
\addto\tenpoint{\normalbaselineskip=1.02\normalbaselineskip\normalbaselines}
\let\=\B
\let\hlop\!
\let\larg\big
\let\largl\bigl
\let\largr\bigr
\let\le\leqslant
\let\ge\geqslant

\let\vt\vartheta

\let\wt\tilde
\let\<\langle
\let\>\rangle
\define\point{\,.\,}

\define\trdeg{\operatorname{tr\,deg}}
\define\Tr{\operatorname{Tr}}
\define\Log{\operatorname{Log}}
\define\alg{\operatorname{alg}}
\define\an{\operatorname{an}}
\define\codim{\operatorname{codim}}
\redefine\det{\operatorname{det}}
\redefine\diag{\operatorname{diag}}
\redefine\mod{\;\operatorname{mod}}
\redefine\pmod#1{\;(\operatorname{mod}#1)}
\define\Chi{\chi^{\vphantom0}}
\define\CC{\Bbb C}
\define\ZZ{\Bbb Z}
\define\QQ{\Bbb Q}

\define\jj{\lower1.8pt\hbox{$\Cal J$}}
\let\kk=k
\define\fH{\frak H}
\define\fK{\frak K}

\define\fR{\frak R}
\define\bA{\boldkey A}
\define\0{\boldkey0}
\define\1{\boldkey1}
\define\ba{\boldkey a}
\define\bb{\boldkey b}
\define\bc{\boldkey c}
\define\bd{\boldkey d}
\define\bn{\boldkey n}

\define\bz{\boldkey z}
\define\bu{\boldkey u}
\define\btau{\boldsymbol\tau}
\define\bpsi{\boldsymbol\psi}
\define\bnu{\boldsymbol\nu}
\define\bdelta{\boldsymbol\delta}
\define\blambda{\boldsymbol\lambda}
\define\bpartial{\boldsymbol\partial}
\define\Eta{{\roman H}}
\define\Aut{\hbox{\sl Aut}}
\define\Hom{\hbox{\sl Hom\/}}
\define\Lie{\hbox{\sl Lie\/}}
\define\Gal{\hbox{\sl Gal}}
\define\Galnabla{\hbox{\sl Gal\/}(\nabla)}
\define\GL{\hbox{\sl GL}}
\define\SP{\hbox{\sl Sp}}
\define\sSP{\hbox{\eightsl Sp}}
\define\Sym{\hbox{\sl Sym}}
\define\Proj{\hbox{\sl Proj\/}}
\define\Der{\hbox{\sl Der\/}}
\define\trans#1{{}^t\!#1}
\define\half{\tfrac12}
\topmatter
\title
On the transcendence degree \\
of the differential field generated \\
by Siegel modular forms
\endtitle
\author
D.~Bertrand\footnotemark\ and W.~Zudilin\footnotemark
\endauthor
\endtopmatter
\footnotetext[1]{Institut de Math\'ematiques de Jussieu
\newline\indent\hbox to2pt{\hss}
{\it E-mail\/}: \tt bertrand\@math.jussieu.fr}
\footnotetext[2]{Moscow Lomonosov State University
\newline\indent\hbox to2pt{\hss}
{\it E-mail\/}: \tt wadim\@ips.ras.ru}
\leftheadtext{D.~Bertrand and W.~Zudilin}
\rightheadtext{Derivatives of modular forms}
\document

\head
\S\,1. Statement of the results
\endhead

It is a classical fact that the elliptic modular function
$\lambda=(\vt_{10}/\vt_{00})^4$ satisfies an algebraic
differential equation of order~$3$ (this goes
back to Jacobi's {\it Fundamenta nova\/}), and
none of lower order (cf.~\cite{Ra},~\cite{M}).
In this paper, we show how these properties generalize to
Siegel modular functions of arbitrary degree.

\medskip

Some notations are necessary before we can state our main results.
Let $g$~ be a positive integer (called indifferently degree,
or genus), let
$\kk$~ be an algebraically closed subfield of~$\CC$, and set

\medskip

\bgroup
\leftskip=8mm \parindent=-\leftskip

$\fH_g$ $=$ Siegel half space of degree~$g$;
the $\QQ$-vector group~$Z_g$ formed by
symmetric matrices of order~$g$ has dimension
$$
n:=\frac{g(g+1)}2,
$$
and $\fH_g$~ is open in~$Z_g(\CC)$;

$\btau=(\tau_{jl})$ a generic point on~$\fH_g$,
so that $\kk(2\pi i\btau)$ can be viewed as the field of
rational functions on~$Z_g/\kk$;

$\Gamma$ $=$ a congruence subgroup of $\SP_{2g}(\ZZ)$
(equivalently, a subgroup of finite index if $g>1$).
We recall that the symplectic group~$\SP_{2g}$
has dimension $\dim\SP_{2g}=2g^2+g$;

$R_w(\Gamma,\kk)$ $=$ $\kk$-vector-space of
$\kk$-rational modular forms of weight~$w$
(a non-negative integer) relative to~$\Gamma$,
i.e\. holomorphic functions~$f$ on~$\fH_g$ which satisfy
$$
f(\gamma\btau)=\det(c\btau +d)^wf(\btau)
\qquad \text{for all} \;
\bigl(\gamma=\bigl(\smallmatrix a&b\\ c&d\endsmallmatrix\bigr),\btau)
\in\Gamma\times\fH_g,
$$
and (if $g=1$) which are holomorphic at the cusps of~$\Gamma$;
such $f$'s admit a Fourier expansion
$$
f(\btau)=\sum_{\bnu\in T^+}a_{\bnu}\exp(2\pi i\Tr(\bnu\btau)),
$$
where $T^+$~is the set of non-negative elements of a suitable
lattice in~$Z_g$, and we require that the coefficients~$a_{\bnu}$
all belong to~$\kk$;

$R:=R(\Gamma,\kk)$ $=$ the graded ring
$\bigoplus_{w\ge0}R_w(\Gamma,\kk)$;
in the fraction field of~$R$, quotients of modular forms
of weights $w_1,w_2$ are called meromorphic modular forms
of weight $w_1-w_2$ (or simply, {\it modular functions\/}
if $w_1=w_2$);

$K:=K(\Gamma,\kk)$ $=$ the field of modular functions;
for $g>1$, the field $K\otimes_{\kk}\CC$ identifies with the
field of meromorphic functions on~$\fH_g$ which are invariant
under the action of~$\Gamma$ (cf.~\cite{S}, \S\,25.4). In fact,
$\Proj(R)$~ is a projective variety over~$\kk$, whose field of
$\kk$-rational functions identifies with~$K$. In particular,
$$
\trdeg(K/\kk)=\frac{g(g+1)}2=n,
\qquad
\trdeg(R/\kk)=n+1,
$$
and we may write $K=\kk(\lambda)$,
where $\lambda:=\{\lambda_1,\dots,\lambda_N\}$~ is a set of
modular functions relative to~$\Gamma$, whose first $n$~elements
$\blambda=\{\lambda_1,\dots,\lambda_n\}$ are algebraically
independent\footnotemark; in particular,

$\partial/\partial\blambda
:=\{\partial/\partial\lambda_1,\dots,\partial/\partial\lambda_n\}$
is a $K$-basis of the space $\Der(K/\kk)$
of $\kk$-derivations of the extension~$K/\kk$;

$\bdelta=\{\delta_{jl},\ 1\le j\le l\le g\}$ where
$$
\delta_{jl}=\frac1{2\pi i}\frac\partial{\partial\tau_{jl}},
\quad 1\le j<l\le g,
\qquad \text{while} \quad
\delta_{jj}=\frac1{\pi i}\frac\partial{\partial\tau_{jj}},
\quad 1\le j\le g.
$$
We sometimes reindex this set of partial
derivations as $\bdelta=\{\delta_1,\dots,\delta_n\}$,
and abbreviate it by $\partial/\pi i\,\partial\btau$;
they form a $\kk(2\pi i\btau)$-basis
of~$\Der(\kk(2\pi i\btau)/\kk)$;

$M:=M(\Gamma,\kk)=\kk\<\lambda_1,\dots,\lambda_N\>$ $=$
the $\bdelta$-differential field generated by~$K$,
i.e\. the field generated over~$\kk$ by the partial derivatives
of all orders in the~$\delta_{jl}$, $1\le j\le l\le g$,
of all the elements of~$K$, or equivalently,
the field generated over~$K$ by the $\bdelta$-derivatives of
all orders of~$\blambda$ (recall that in characteristic~$0$,
the derivatives of an element~$x$ algebraic over a differential
field~$K$ lie in~$K(x)$).

\egroup

\footnotetext[3]{As a rule, we reserve boldface Greek
letters to $n$-tuples (of numbers, functions or
derivatives).}

\medskip

With these notations in mind, we can state:

\proclaim{Theorem 1}
The $\bdelta$-differential field $M=\kk\<\lambda\>$
generated by the field of modular functions
$K=\kk(\lambda)$~is a finite extension of the field
generated over~$K$ by the $\bdelta$-partial derivatives
of $\lambda_1,\dots,\lambda_n$ of order~$\le2$\rom,
and has transcendence degree
$$
\trdeg(M/\kk)=\dim\SP_{2g}=2g^2+g
$$
over~$\kk$. Furthermore\rom, $M$ and $\CC(\btau)$
are linearly disjoint over~$\kk$\rom, hence
$$
\trdeg(M(2\pi i\btau)/\kk)
=\dim(\SP_{2g}\times Z_g)
=\frac12g(5g+3)
$$
\rom(and $\pi$~is transcendental over $M(2\pi i\btau)$ if
$\kk=\QQ^{\alg}$\rom)\rom, while all modular forms
in~$R$ are algebraic over~$M$.
\endproclaim

The statement concerning~$M$ and~$\CC(\btau)$ is clear,
since $M$ embeds in the fraction field of the ring of
convergent Puiseux series in $\exp(2\pi i\tau_{jl})$,
$1\le j\le l\le g$, with coefficients in~$\kk$,
which is linearly disjoint from~$\CC(\btau)$ over~$\kk$.
Thus, the second formula is an immediate corollary of the first one.
Now, both are easily seen to be equivalent,
and our strategy will consist in proving
the second formula, i.e\. in studying the
$\partial/\pi i\,\partial\btau$-differential
extension $M(2\pi i\btau)$ of~$\kk$.

Since an algebraic extension of a differential field of
characteristic~$0$ is automatically a differential extension,
we may assume, in order to prove Theorem~1, that $\Gamma$~ is
contained in a principal congruence subgroup of level
at least~$3$, so that $\Gamma\backslash\fH_g$~ is a complex
manifold, whose natural image in $\Proj(R(\Gamma,\kk))$~ is
the set of complex points of a smooth quasiprojective variety
$S(\Gamma)/\kk$ (cf.~\cite{MF}, p.~190, and~\cite{P}, \S\,2).
For instance, we may take for~$\Gamma$ the theta-group of level~$(4,8)$
$$
\Gamma_{4,8}
=\bigl\{\gamma=\bigl(\smallmatrix a&b\\c&d\endsmallmatrix\bigr)
\in\SP_{2g}(\ZZ),\ \gamma\equiv\1_g\pmod4,\
\diag(a\trans b)\equiv\diag(c\trans d)\equiv0\pmod8\bigr\}.
$$
By a well-known result of Igusa (\cite{I}, pp.~178,~190 and~224),
the corresponding ring $R(\Gamma_{4,8},\kk)$~ is
the integral closure of the ring
$\kk[\vt_{\ba}\vt_{\bb},\
\ba\in(\ZZ/2\ZZ)^{2g},\ \bb\in(\ZZ/2\ZZ)^{2g}]$,
where
$$
\vt_{\ba}(\btau)
=\vt_{(\ba',\ba'')}(\btau)
=\sum_{\bn\in\ZZ^g}
\exp\bigl(\pi i\bigl(\trans{\bigl(\bn+\half\ba'\bigr)}
\btau\bigl(\bn+\half\ba'\bigr)
+\trans{\bigl(\bn+\half\ba'\bigr)}\ba''\bigr)\bigr)
$$
denotes the `thetanull' modular form attached to the 2-characteristic
$\ba=(\ba',\ba'')\in\allowmathbreak(\ZZ/2\ZZ)^{2g}$.
In particular, we may here take
$\lambda=\{\vt_{\ba}/\vt_{\0},\ \ba\in(\ZZ/2\ZZ)^{2g}\}$
as a set of generators of~$K(\Gamma_{4,8},\kk)$.
In this case, Theorem~1 can be given a more precise form, as follows.

\proclaim{Theorem 2}
The $\bdelta$-derivatives of order~$\le2$ of the modular functions
$\{\vt_{\ba}/\vt_{\0},\ \ba\in\allowmathbreak(\ZZ/2\ZZ)^{2g}\}$
generate over~$\kk$ a $\bdelta$-stable field $M(\Gamma_{4,8},\kk)$
of transcendence degree $2g^2+g$ over~$\kk$\rom,
over which $\vt_{\0}$~ is algebraic.
\endproclaim

In fact, the thetanulls themselves satisfy a system of
partial differential equations of the {\it second\/} (rather
than third) order, which has been explicited in~\cite{Z},
Theorem~1, and is discussed, amongst other differential
properties, in~\S\,5 of the paper. As shown by Ohyama~\cite{O},
these differential equations take a simpler form in the case
$g=2$. In~\S\,6, we derive from these results and from
Theorem~2 ten explicit functions making up a transcendence
basis for the corresponding field~$M$, as well as an
algebraic presentation of the $\bdelta$-stable {\it ring\/}
generated over~$\QQ$ by the ten non-zero degree~$2$ thetanulls
$\vt_{\ba}$, and their thirty first order partial logarithmic derivatives
$\psi_{\ba,l}=\delta_l\vt_{\ba}/\vt_{\ba}$. See Theorem~3
of~\S\,6 for a precise statement.

\medskip

The following diagram and `legend' illustrate the proof of Theorem~1.
The notation~$P$ is defined in~\S\,3, Lemma~2, and $\asymp$
means equality for the algebraic closures~$\bullet^{\alg}$.
We show that $\wt M^{\alg}=M^{\alg}$ in~\S\,4.
But as hinted above, the crux of the proof consists
in translating the problem in terms of a `field of periods'
$F:=\wt M(2\pi i\btau,2\pi i)$ and its compositum~$\Phi$ with~$\CC$,
hence in terms of {\it linear\/} differential equations,
as described in~\S\,2. The fact that the (algebraic) differential
equations satisfied by classical modular forms are governed
by the (linear) Gauss--Manin connection was made
crystal clear in~\cite{Ka}, Appendix~1, and we are here
merely extending this view-point to forms of higher degrees.

\nopagebreak
\bgroup
\input picture.sty
\unitlength=9mm
\begin{picture}(0,5.25)(-7,-.1)
\put(-5.55,4.7){\makebox(0,0){$\Phi=\CC(\lambda,\Pi(\lambda))$}}
\put(-3.2,4.7){\makebox(0,-.18){$=_{[\otimes_{\kk(2\pi i)}\CC]}$}}
\put(-3.2,2.2){\makebox(0,-.18){$\kk(2\pi i\btau)$}}
\put(-2.85,2.6){\vector(1,2){.8}}
\put(-3.55,2.6){\vector(-1,2){.8}}
\put(-.95,.15){\vector(-3,2){2}}
\put(-.3,0){\makebox(0,0){$\kk$}}
\put(-.3,.4){\vector(0,1){.6}}
\put(-.3,1.4){\makebox(0,0){$\kk(P)$}}
\put(-.3,1.8){\vector(0,1){.6}}
\put(-.3,3){\makebox(0,0){$\wt M:=\kk(\lambda,P,\bdelta P)$}}
\put(-.3,3.6){\vector(0,1){.6}}
\put(-.6,4.7){\makebox(0,0){$F:=\wt M(2\pi i,2\pi i\btau)$}}
\put(3.85,4.7){\makebox(0,0){$M(2\pi i\btau)$}}
\put(2.1,4.7){\makebox(0,-.18){$\asymp_{[\otimes_{\kk}\kk(2\pi i)]}$}}
\put(4.8,3.55){\vector(-3,2){1}}
\put(4.6,4){\makebox(0,0){$\ssize Z_g$}}
\put(4.7,3){\makebox(0,0){$M=\kk\<\lambda\>$}}
\put(2.5,3){\makebox(0,0){$\asymp$}}
\put(4.7,1.85){\vector(0,1){.7}}
\put(4.7,1.4){\makebox(0,0){$K=\kk(\lambda)$}}
\put(.35,.15){\vector(3,1){3}}
\end{picture}
\egroup

\vskip-3mm
$$
\align
\Gal_{\partial/\partial\blambda}(\Phi/K\otimes_{\kk}\CC)=\SP_{2g}
&\implies \trdeg(M/K)=\codim_{\sSP_{2g}}(Z_g)
\\
&\implies \trdeg(M/\kk)=\dim\SP_{2g}
\endalign
$$

Further comments on the proof of Theorem~1, and on its
possible generalizations and applications, are given
in Remarks~3 and~4 at the end of~\S\,4. See also~\S\,3,
Remark~2 for a direct relation between periods and
derivatives of modular forms {\it via\/} modular tensors.

\head
\S\,2. Picard--Fuchs and Picard--Vessiot theories
\endhead

We first recall some well-known facts (see~\cite{D2}, \S\,2)
about algebraic families of abelian varieties. Let
$S/\kk$~ be a smooth algebraic variety over $\kk\subset\CC$,
with generic point~$\sigma$ and field of rational functions
$L=\kk(S)=\kk(\sigma)$, and let $f\:\bA\to S$~ be
a principally polarized abelian scheme over~$S$
of relative dimension~$g$, with generic fiber
$\bA_\sigma:=A$ $=$ an abelian variety over~$L$.
We write $f^{\an}$ for the analytic map deduced
from~$f$ after extension of scalars to~$\CC$.
Let $H:=H^1_{dR}(A/L)$ be the $2g$-dimensional $L$-vector space
formed by the cohomology classes of $L$-rational
differential forms of the second kind on~$A/L$.
The Gauss--Manin connection attached to~$f$ equips~$H$
with an integrable connection $\nabla\:H\to H\otimes
\Omega^1_{L/k}$. Choose a base point $s_0\in S^{\an}$, and
let $L_0$ be the field of meromorphic functions over a
small neighbourhood $U(s_0)$ of~$s_0$. The restriction
to~$U(s_0)$ of the local system formed by the relative Betti
cohomology $\bold R^1f^{\an}_*\QQ$  generates over~$\CC$ the full
space of horizontal vectors of the extension of~$\nabla$
to~$H\otimes L_0$, and provides
a representation~$\rho$ of the fundamental group $\pi_1(S^{\an},s_0)$
on the $2g$-dimensional $\QQ$-vector space
$H_B=H^1_B(\bA_{s_0}^{\an},\QQ)$, which preserves the
$2\pi i\QQ$-valued symplectic form~$\psi_B$
induced on~$H_B$ by the principal polarisation
on~$\bA_{s_0}$. The Gauss--Manin connection admits an
extension with logarithmic singularities over a suitable
compactification of~$S$, so that $\nabla$~ is fuchsian
and the Zariski closure of $\rho(\pi_1)$ in
$\Aut_{\psi_B}(H_B\otimes\CC)$ is isomorphic to the {\it
differential Galois group\/}
$\Galnabla$ of the connection~$\nabla$.

Now, choose a basis $\omega_1,\dots,\omega_g$,
$\eta_1,\dots,\eta_g$ of~$H$ over~$L$,
and a basis $c_1,\dots,c_{2g}$ of the relative Betti homology
of $\bA^{\an}$ above~$U(s_0)$.
A basis of horizontal vectors of~$\nabla$
is then represented by the inverse of the
fundamental matrix of periods and quasi-periods
$$
\Pi(\sigma)
=\pmatrix \int_{c_i}\omega_j \\ \int_{c_i}\eta_j \endpmatrix
=\pmatrix \Omega_1(\sigma) & \Omega_2(\sigma) \\
\Eta_1(\sigma) & \Eta_2(\sigma) \endpmatrix,
$$
whose coefficients extend over~$S^{\an}$ to multivalued
meromorphic functions in~$\sigma$. By definition, they generate over
$L\otimes_{\kk}\CC=\CC(\sigma)$ a Picard--Vessiot extension
$$\Phi=\CC(\sigma,\Pi(\sigma)).$$
In particular, $\Phi$~ is stable under the partial derivatives
$\bpartial=\{\partial_1,\dots,\partial_{\dim S}\}$
given by a basis over~$L$ of the dual $\Der(L/\kk)$ of
$\Omega^1_{L/\kk}$, and by the main theorem of Picard--Vessiot theory,
$\trdeg(\Phi/\CC(\sigma))=\dim\Galnabla$.

Assume now that the family~$f$ is `as generic as possible',
or more precisely (cf.~\cite{D1}, \S\,4.4.14.1) that the
induced morphism~$\phi_f$ from~$S$ to one of the components
of the moduli space of principally polarized abelian
varieties is {\it dominant}. We then have:

\proclaim{Lemma 1 \rm(\cite{D1}, Lemme~4.4.16)}
The image of~$\pi_1(S,s_0)$ under~$\rho$ has finite index
in the group $\Aut_{\psi_B}(H_B)$ of symplectic automorphisms
of~$H_B$.
\endproclaim

The genericity hypothesis holds tautologically
when $S$~ is the moduli scheme of principally
polarized abelian varieties with a level
structure of order~$\ge3$, hence, in the notations of~\S\,1,
when $S=S(\Gamma)$, endowed with the
corresponding universal abelian scheme~$f$
(cf.~\cite{MF}, Appendix~7, A--B).
In this case, $L$~ is isomorphic to $K=K(\Gamma,\kk)=\kk(\lambda)$,
and we shall write indifferently $\lambda$ for~$\sigma$, e.g.
$$
\Pi(\sigma)=\Pi(\lambda),
\qquad
\bpartial=\frac\partial{\partial\blambda}
:=\{\partial/\partial\lambda_1,\dots,\partial/\partial\lambda_n\},
$$
while the multivalued map $\lambda\mapsto\Pi(\lambda)$ lifts
(after a choice of a base point~$\btau_0$ above~$s_0$)
to a meromorphic map
$$
\btau\in\fH_g\mapsto\wt\Pi(\btau)
:=\Pi(\lambda(\btau))
=\pmatrix \wt\Omega_1(\btau) & \wt\Omega_2(\btau) \\
\wt\Eta_1(\btau) & \wt\Eta_2(\btau) \endpmatrix
$$
on the universal covering manifold~$\fH_g$
of $S^{\an}\simeq\Gamma\backslash\fH_g$.
Lemma~1 then implies that
$\Gal_{\partial/\partial\blambda}(\Phi/K\otimes\CC)
\hlop=\hlop\Galnabla
\hlop=\hlop\Aut_{\psi_B}(H_B\otimes\CC)
\hlop\simeq\hlop\SP_{2g}(\CC)$
so that
$\trdeg(\Phi/\CC(\lambda))\hlop=\allowmathbreak\dim\SP_{2g}$,
hence

\proclaim{Proposition 1}
Assume that $S=S(\Gamma)$\rom, where $\Gamma$~ is any congruence
subgroup of~$\SP_{2g}(\ZZ)$. Then\rom,
$$
\trdeg\bigl(\CC(\lambda,\Pi(\lambda))/\CC\bigr)
=\dim(\SP_{2g}\times Z_g).
$$
\endproclaim

\demo{Proof}
$\trdeg(\Phi/\CC)
=\trdeg(\Phi/\CC(\lambda))+\trdeg(\CC(\lambda)/\CC)
=\dim\SP_{2g}+n$,
which is $\dim(\SP_{2g}\times Z_g)$.
\enddemo

\remark{Remark \rom1}
In the sequel, the intermediate fields
$$
\Phi_{\kk}:=\kk(\lambda,\Pi(\lambda)),
\qquad
\Psi_{\kk}:=\kk\largl(\lambda,\frac1{2\pi i}\Omega_1(\lambda),
\frac1{2\pi i}\Eta_1(\lambda)\largr)
$$
will be used. Since $\nabla_\partial$ acts on the
$K$-vector space~$H$ for any  $\partial\in\Der(K/\kk)$,
they are still $\partial/\partial\lambda$-differential
extensions of~$K$, but even $\Phi_{\kk}$ is in general
{\it not\/} a Picard--Vessiot extension.
Indeed, the generalized Riemann relations
(i.e\. the reciprocity law for differentials of
the second kind on~$A$; cf.~\cite{B}, pp.~37--38,
or~\cite{D2}, Proposition~1.5, in connection with the
form~$\psi_{dR}$ introduced below) show that the
$\partial/\partial\blambda$-constant~$\pi$ lies
in~$\Phi_{\kk}$, although not in~$K$ if $\kk=\QQ^{\alg}$.
In a sense, it is the field $\Psi_{\kk}(2\pi i\btau)$ which
gives the required Picard--Vessiot extension,
but it is not a field of periods in the usual sense:
it will in general contain neither~$\btau$, nor~$2\pi i$\,!
Also, note that these Riemann relations
(written in the standard bases) immediately imply that
$\trdeg(\Phi/\CC(\sigma))\le2g^2+g$
for any family $\bA\to S$: once $\Omega_1$~ is chosen
(at most $g^2$~degrees of freedom), they leave at most
$g(g+1)/2$ degrees of freedom for the entries of~$\Omega_2$,
at most $g(g+1)/2$ for those of~$\Eta_1$, and none for~$\Eta_2$.
\endremark

\head
\S\,3. Periods of the first kind and their derivatives
\endhead

Let $f\:\bA\to S$ be a principally polarized abelian
scheme~$S$ as in~\S\,2. In order to prepare for a modular study of
the period matrix $\Pi(\sigma)$, we shall need some
well-known infinitesimal properties of the structural
morphism~$\phi_f$. Recall that in parallel with~$\psi_B$,
the polarisation on~$\bA_\sigma$ provides a
non-degenerate antisymmetric form~$\psi_{dR}$ on
$H:=H^1_{dR}(A/L)$ with values in~$L$, which admits as
maximal isotropic subspace the $L$-vector space
$\Omega:=\Omega^1_{A/L}$ of~$H$ formed by the (cohomology
classes of) differentials of the first kind on~$A/L$,
so that $H/\Omega$ is canonically isomorphic to the $L$-dual
$\Omega^*\simeq\Lie(A/L)$ of~$\Omega$.

When writing $\Pi(\sigma)$,
we may choose a $\psi_B$-symplectic basis%
\footnote[4]{We order this basis in such a way that
$\btau=\Omega_1^{-1}\Omega_2\in\fH_g$; then,
$\Gamma$~ acts on the $(g\times 2g)$-matrix
$(\wt\Omega_2(\btau) \; \wt\Omega_1(\btau))
=\wt\Omega_1(\btau)\,(\btau \; \1_g)$
by the {\it transpose\/} of its standard representation.}
of the Betti homology for $\{c_1,\dots,c_{2g}\}$,
and a basis of~$\Omega$ over~$L$ for $\omega_1,\dots,\omega_g$.
Then, $\Omega_1(\sigma)$ is invertible, and
$\btau(\sigma):=\Omega_1(\sigma)^{-1}\Omega_2(\sigma)$
is one of the points in~$\fH_g$ parametrizing the principally
polarized abelian variety~$\bA_\sigma$. As for the last
$g$~rows of~$\Pi(\sigma)$, we describe them with the help of
the Kodaira--Spencer map
$$
\kappa_f\:\Der(L/k)\to\Hom(\Omega,H/\Omega)
\simeq(\Omega^*)^{\otimes 2}\:
\partial\mapsto\{\omega\mapsto\nabla_{\partial}(\omega)\mod\Omega\}
$$
attached to~$f$ (cf.~\cite{L}, p.~157).
Note that the image of~$\kappa_f$ is contained in
$\Sym^2(\Omega^*)$, because $\psi_{dR}$ is horizontal for~$\nabla$.

Assume now that $f$ is `sufficiently generic' in
the sense of~\S\,2. Since $\kappa_f$ represents the tangent map to
the morphism~$\phi_f$ at the generic point~$\sigma$ of~$S$,
its rank then coincides with $g(g+1)/2$ and $\kappa_f$ maps
onto $\Sym^2(\Omega^*)$ (cf.~\cite{K}, p.~169, and \cite{P}, p.~255).
Therefore, there exists a partial derivation $\partial_0\in\Der(L/\kk)$
whose image under~$\kappa_f$ is an isomorphism from~$\Omega$
to~$H/\Omega$. In particular, the differential forms of the second kind
$\{\eta_l=\nabla_{\partial_0}(\omega_l),\ 1\le l\le g\}$
lift a basis of $H/\Omega$ over~$L$. Since
$\int_c\nabla_\partial(\omega)=\partial(\int_c\omega)$
for any $\partial\in\Der(L/\kk)$, $\omega\in H$, $c\in H_B$,
we then infer that $\Eta_1(\sigma)$ (resp\.
$\Eta_2(\sigma)$) are the $\partial_0$-derivatives of
$\Omega_1(\sigma)$ (resp\. $\Omega_2(\sigma)$). In other words,
$\Pi(\sigma)$ takes in such bases the shape
$$
\Pi(\sigma)
=\pmatrix
\Omega_1(\sigma) & \Omega_1(\sigma)\btau(\sigma) \\
\partial_0\Omega_1(\sigma) & \partial_0(\Omega_1(\sigma)\btau(\sigma))
\endpmatrix.
$$
This certainly implies that the
$\bpartial$-stable field $\Phi=\CC(\sigma,\Pi(\sigma))$~ is
generated over $\CC(\sigma)$ by $\Omega_1(\sigma),\btau(\sigma)$,
$\bpartial\Omega_1(\sigma),\bpartial\btau(\sigma)$
(recall that $\bpartial$ stands for a basis of~$\Der(L/\kk)$),
and more precisely, that the $\bpartial$-stable
field~$\Phi_{\kk}$ of Remark~1 can be written as
$$
\Phi_{\kk}:=\kk(\sigma,\Pi(\sigma))
=\kk(\sigma,\Omega_1(\sigma),
\btau(\sigma),\bpartial\Omega_1(\sigma),
\bpartial\btau(\sigma)).
$$
Similarly, each of the columns of the fundamental matrix
of solutions $\Pi(\sigma)$ generates over~$L$ a $\bpartial$-stable
field, and by the same argument, the $\bpartial$-stable field
$\Psi_{\kk}$ of Remark~1 reads
$$
\Psi_{\kk}
:=\kk\largl(\sigma,\frac1{2\pi i}\Omega_1(\sigma),
\frac1{2\pi i}\Eta_1(\sigma)\largr)
=\kk\largl(\sigma,\frac1{2\pi i}\Omega_1(\sigma),
\bpartial\largl(\frac1{2\pi i}\Omega_1(\sigma)\largr)\largr).
$$

\medskip

In order to compare $\Phi_{\kk}$ with the (still to be defined)
fields~$F$ and $\wt M$ of our diagram, we restrict from now
on to the modular situation of Proposition~1, with
$S=S(\Gamma)$, $\sigma=\lambda$, $L=K$,
$\wt\Omega_1(\btau)=\Omega_1(\lambda(\btau))$, etc.
It will be useful (though not strictly necessary, cf\. Remark~2 below)
to choose an explicit basis of the $K$-vector space~$\Omega$.
We appeal to Shimura's differentials (\cite{S}, \S\,30)
for such a specification, and to fix notations,
henceforth assume that $\Gamma =\Gamma_{4,8}$,
so that $S(\Gamma)$ is the moduli scheme of principally
polarized abelian varieties with $(4,8)$-level
structure (cf.~\cite{MF}, pp.~193--195).

Consider the full set of abelian theta functions with two-characteristics
$\{\vt_{\ba}(\bz,\btau),\allowmathbreak\ba\in(\ZZ/2\ZZ)^{2g}\}$,
as given in~\S\,5 below. Then,
$\{\bz\mapsto\vt_{\ba}(2\bz,\btau),\ \ba\in(\ZZ/2\ZZ)^{2g}\}$
defines a projective embedding of
$\CC^g/(\ZZ^g\oplus\btau\ZZ^g)$ (cf.~\cite{I}, pp.~169 and~171),
whose image can be identified with the generic fiber
$\bA_{\lambda(\btau)}$ of the corresponding universal
family. Therefore (cf.~\cite{S}, Lemma~30.2), we can
choose
$\vt_{\0}=\vt_{(\0,\0)}$ and $g$~ odd theta functions
$\vt_1,\dots,\vt_g$ amongst them such that the
jacobian matrix at $\bz=\0$ of the map
$$
\bz\mapsto\largl(\frac{\vt_j}{\vt_{\0}}
\largl(\frac{\bz}{2\pi i},\btau\largr),\
j=1,\dots,g\largr)
$$
is invertible. Because of parities,
this jacobian reads $\trans{P(\btau)}$, where
$$
P(\btau)
=\frac1{2\pi i\,\vt_{\0}(\0,\btau)}
\cdot\pmatrix
\dfrac{\partial\vt_1}{\partial z_1}(\0,\btau) & \dots &
\dfrac{\partial\vt_1}{\partial z_g}(\0,\btau)
\\ \noalign{\vskip4pt}
\hdotsfor3
\\ \noalign{\vskip4pt}
\dfrac{\partial\vt_g}{\partial z_1}(\0,\btau) & \dots &
\dfrac{\partial\vt_g}{\partial z_g}(\0,\btau)
\endpmatrix.
$$

\proclaim{Lemma 2 \rm(\cite{S}, Theorem~30.3 and Formula~(30.2f))}
Let $\btau\in\fH_g$\rom, and let $\lambda=\lambda(\btau)$.
The relation
$$
(\omega_1,\dots,\omega_g)
=(dz_1,\dots,dz_g)\,2\pi i\,\trans{P(\btau)}
$$
defines a basis of differential forms of the first kind on
the abelian variety $\bA_{\lambda(\btau)}$ which are
rational over the field $K=\kk(\lambda)$\rom, and which admit
$2\pi iP(\btau)\,(\1_g \; \btau)$ as period matrix.
In particular\rom, $(1/2\pi i)\,\wt\Omega_1(\btau)=P(\btau)$
in such a basis.
\endproclaim

Recall the notations at the beginning of~\S\,1,
and set further

\medskip

\bgroup
\leftskip=8mm \parindent=-\leftskip

$\xi$ $=$ the map from $\Phi=\CC(\lambda,\Pi(\lambda))$
to the field of meromorphic functions on~$\fH_g$ which
lifts an element $f(\lambda)\in\Phi$ to
$$
(\xi(f))(\btau)=\wt f(\btau):=f(\lambda(\btau));
$$
relatively to the bases
$\{\partial/\partial\blambda\}$ and
$\{\bdelta=\partial/\pi i\,\partial\btau\}$,
the differential of the covering
map~$\lambda$ corresponding to~$\xi$ is given by
$(\delta_1,\dots,\delta_n)
=(\partial/\partial\lambda_1,\dots,\allowmathbreak
\partial/\partial\lambda_n)\trans{W(\btau)}$,
where $W(\btau)$ denotes the invertible matrix
$$
W(\btau)
:=\biggl(\frac{\partial\blambda}{\pi i\,\partial\btau}\biggr)(\btau)
=\pmatrix
\delta_1\lambda_1 & \dots & \delta_1\lambda_n \\
\hdotsfor3 \\
\delta_n\lambda_1 & \dots & \delta_n\lambda_n
\endpmatrix;
$$

$\wt M:=\kk\bigl(\lambda(\btau),P(\btau),
\bdelta P(\btau),\bdelta\lambda(\btau)\bigr)
=\kk\bigl(\lambda(\btau),P(\btau),\bdelta P(\btau),W(\btau)\bigr)$
(since $\lambda_{n+1},\dots,\allowmathbreak\lambda_N$
are algebraic over the field~$K$).
A simpler description of~$\wt M$
will presently be given, cf\. Remark~2 below;

$F:=\wt M(2\pi i\btau,2\pi i)$.

\egroup

\medskip

The field~$F$ being thus defined, we can now relate it
to the fields of periods~$\Phi$ of~\S\,1, as follows.

\proclaim{Proposition 2}
Assume that $S=S(\Gamma)$\rom, with $\Gamma=\Gamma_{4,8}$.
Then\rom, $\xi$ induces an isomorphism from
$\Phi_{\kk}=\kk(\lambda,\Pi(\lambda))$ onto~$F$.
Moreover\rom, both fields $F$ and $\wt M$
are stable under the partial derivatives
$\bdelta=\partial/\pi i\,\partial\btau$.
\endproclaim

\demo{Proof}
We already know that
$$
\xi(\Phi_{\kk})
=\kk\largl(\lambda(\btau),2\pi i\,P(\btau),\btau,
2\pi i\,\xi\largl(\frac\partial{\partial\blambda}P(\lambda)\largr),
\xi\largl(\frac\partial{\partial\blambda}\btau(\lambda)\largr)\largr),
$$
and that it contains~$2\pi i$. Now, for any $f\in\Phi$,
$$
\bigl(\delta_1(\xi(f)),\dots,\delta_n(\xi(f))\bigr)
=\largl(\xi\largl(\frac{\partial f}{\partial\lambda_1}\largr),\dots,
\xi\largl(\frac{\partial f}{\partial\lambda_n}\largr)\largr)
\trans{W(\btau)};
$$
for instance,
$\pi i\xi((\partial/\partial\blambda)\btau(\lambda))$
can be written as the inverse of the matrix~$W(\btau)$,
and they have the same field of definition.
Therefore, $\xi(\Phi_{\kk})$~ is generated over~$\kk$
by $2\pi i$, $2\pi i\btau$, $\lambda(\btau)$, $P(\btau)$,
$\bdelta P(\btau)$, and $W(\btau)$, which precisely
constitute a set of generators for~$F/\kk$. As for the second
part of Proposition~2, note that $\Phi_{\kk}$~ is
$\partial/\partial\blambda$-stable, so that
$\xi(\Phi_{\kk})=F$~ is
$(\partial/\pi i\,\partial\btau)\trans{W(\btau)^{-1}}$-stable,
hence $\partial/\pi i\,\partial\btau$-stable as well,
since the coefficients of $\trans{W(\btau)^{-1}}$ lie in~$F$.
Similarly, we know that the subfield~$\Psi_{\kk}$ of~$\Phi_{\kk}$
generated over~$K$ by the first $g$~columns of
$(1/2\pi i)\,\Pi(\lambda)$~ is stable under
$\partial/\partial\blambda$, and the same argument implies that
$\wt M=\xi(\Psi_{\kk})\bigl(W(\btau)\bigr)$ too is
$\partial/\pi i\,\partial\btau$-stable.
\enddemo

\remark{Remark \rom2}
In view of~\cite{S}, Formula~30.2d (or more simply,
of the monodromy action on periods, cf\. Footnote~$^{(4)}$),
$\trans P$~ is a modular tensor relative to the standard
representation $\rho$~of $\GL_g(\CC)$, i.e.:
$$
\trans{P}(\gamma\btau)
=(c\btau+d)\trans{P}(\btau)
\qquad \text{for all} \;
\gamma=\bigl(\smallmatrix a&b\\c&d\endsmallmatrix\bigr)\in\Gamma.
$$
On the other hand, it is a well-known fact that the
$\bdelta$-derivatives of any modular function~$\lambda_0$
can be arranged into a (meromorphic) $\Sym^2\rho$-modular tensor;
in other words (cf\. Formula~\thetag{4} of~\S\,5 below):
$$
(\bdelta\lambda_0)(\gamma\btau)
=(c\btau+d)\cdot(\bdelta\lambda_0)(\btau)\cdot\trans{(c\btau+d)}
\qquad \text{for all} \;
\gamma=\bigl(\smallmatrix a&b\\c&d\endsmallmatrix\bigr)\in\Gamma.
$$
Consequently, the matrix valued functions
$\trans{P^{-1}}\bdelta\lambda_jP^{-1}$, $j=1,\dots,n$,
are invariant under~$\Gamma$.
Looking at their Fourier expansions, we deduce that
their entries all belong to the field
$K=K(\Gamma,\kk)=\kk(\lambda)$. Therefore, the last
set of generators $\bdelta\lambda(\btau)$ (equivalently: the
entries of~$W(\btau)$) occuring in the definition of~$\wt M$
already lies in the field~$K(P(\btau))$ generated by the
first two ones, i.e\.
$\kk(\lambda,\bdelta\lambda)\subset\kk(\lambda,P)$,
and in parallel with $F=\xi(\Phi_{\kk})$, we finally obtain
$$
\wt M=\kk(\lambda(\btau),P(\btau),\bdelta P(\btau))
=\kk\largl(\lambda(\btau),\frac1{2\pi i}\wt\Omega_1(\btau),
\frac1{2\pi i}\wt\Eta_1(\btau)\largr)
=\xi(\Psi_{\kk})
$$
as a simpler expression for the field~$\wt M$
of the diagram of~\S\,1.
\endremark

In fact, this type of argument can be reversed
(using the fact that $\bdelta\lambda_1,\dots,\bdelta\lambda_n$
form a {\it basis\/} of the $K$-vector-space of meromorphic
$\Sym^2\rho$-modular tensors), and implies that all
binomials in the entries of~$P(\btau)$ belong to the
field~$K(\bdelta\lambda)$. In particular,
$k(\lambda,P)\subset(\kk(\lambda,\bdelta\lambda))^{\alg}$.
We shall give a precise version of this statement
in Proposition~4 of~\S\,5, but notice that the above
proofs and results extend to the study of
$(1/2\pi i)\,\wt\Omega_1(\btau)$,
$(1/2\pi i)\,\wt\Eta_1(\btau)$
for any congruence subgroup~$\Gamma$ of~$\SP_{2g}(\ZZ)$.

\head
\S\,4. Proof of Theorems~1 and ~2
\endhead

We now complete the proof of Theorems~1 and~2 in four steps.
Since the algebraic closure of $K(\Gamma,\kk)$,
hence of $M(\Gamma,\kk)$, is independent of the group~$\Gamma$,
Theorem~1 is a corollary of Theorem~2 and for simplicity,
we throughout assume that $\Gamma=\Gamma_{4,8}$,
although this hypothesis may truly be needed only in the last step:
the first three ones are valid for all $\Gamma$'s,
and suffice for the proof of Theorem~1.

For each $m=0,1,\dots,\infty$, we denote by $K^{(m)}$
the field generated over~$\kk$ by the partial derivatives
of order $\le m$ with respect to $\bdelta=\partial/\pi i\,\partial\btau$
of all the elements of~$K$. Thus,
$K^{(0)}=K\subset K^{(1)}=K(\bdelta\lambda)\subset\dots
\subset K^{(\infty)}=M$.
All these fields, as well as $\wt M$, $F$ and the ring~$R$,
are contained in the field of meromorphic functions on~$\fH_g$.

\proclaim{Step 1}
$R\subset(K^{(1)})^{\alg}$ \rom(hence
$R\subset M^{\alg}$\rom).
\endproclaim

It suffices to prove that $K^{(1)}$ contains
a non-zero meromorphic modular form~$\Delta$
of positive weight~$w$, since for any modular form~$f$
of weight~$w(f)$, $f^w/\Delta^{w(f)}$ will then lie in~$K$.
(Incidentally, this further implies that the logarithmic derivatives
$\bdelta f/f$ of~$f$ belong to~$K^{(2)}$.)
An explicit choice for~$\Delta$ is given in Lemma~4 of~\S\,5,
but here is a general construction, along the lines of Remark~2 above.
Since the $\bdelta$-derivatives of any modular function~$\lambda_0$
make up a (meromorphic) $\Sym^2\rho$-modular tensor $\bdelta\lambda_0$
with respect to the standard representation~$\rho$ of~$GL_g(\CC)$,
and since the representation
$\Lambda^{g(g+1)/2}(\Sym^2\rho)$ of~$\GL_g(\CC)$ is
isomorphic to $(\Lambda^g\rho)^{\otimes(g+1)}$,
the exterior product $\Delta(\btau)$ of
$\bdelta\lambda_1,\dots,\bdelta\lambda_n$
is a non-zero, meromorphic modular form of weight~$g+1$.
Now, $\Delta(\btau)$~ is given in coordinates by the determinant
$\det(\partial\blambda/\pi i\,\partial\btau)$
of the matrix $W(\btau)$ introduced in~\S\,3,
which clearly lies in $\kk(\lambda,\bdelta(\blambda))=K^{(1)}$.

\proclaim{Step 2}
$K^{(2)}\subset M\subset\wt M\subset(K^{(2)})^{\alg}$
\rom(hence $M$~ is finite over~$K^{(2)}$\rom).
\endproclaim

The first inclusion is trivial. Since $M$~ is
the $\delta$-differential field generated by~$K$,
the second one follows from the stability of $\wt M\supset K$
under differentiation (cf\. Proposition~2).
To check the last one, we appeal to Proposition~4~(i) below,
according to which
$$
\largl(\frac1{2\pi i\,\vt_{\0}}\frac{\partial\vt_{\ba}}{\partial z_j}
\larg|_{\bz=\0}\largr)^2
=\frac1{2^{g-1}}\sum_{\bb\in(\ZZ/2\ZZ)^{2g}}(-1)^{\trans{\ba'}\bb''}
\largl(\frac{\vt_{\ba+\bb}}{\vt_{\0}}\largr)^2
\largl(\frac{\vt_{\bb}}{\vt_{\0}}\largr)^2
\frac{\delta_{jj}(\vt_{\ba+\bb}/\vt_{\0})}
{\vt_{\ba+\bb}/\vt_{\0}}
$$
for any odd theta function~$\vt_{\ba}$.
Since quotient of thetanulls are modular functions for
$\Gamma_{4,8}$, this implies that the entries of $P(\btau)$
are algebraic over~$K^{(1)}$
(see the end of Remark~2 for an implicit argument).
Both $P$ and $\bdelta P$ (and of course $\bdelta\lambda$)
are then algebraic over~$K^{(2)}$, and $\wt M$, hence $M$,
is a finite extension of~$K^{(2)}$.

\proclaim{Step 3}
$\trdeg(M(2\pi i\btau)/\kk)=\dim(\SP_{2g}\times Z_g)$.
\endproclaim

From the explicit expression of the Fourier expansions
of its generators, we infer that $\wt M$ embeds in the fraction field
of the ring of convergent Puiseux series in $\exp(2\pi i\tau_{jl})$,
$1\le j\le l\le g$, with coefficients in~$\kk$. Since this
field is linearly disjoint from $\CC(\btau)$ over~$\kk$, and
since $M$ and $\wt M$ have the same algebraic closure by Step~2,
we deduce from Proposition~2 that
$$
\align
\trdeg(M(2\pi i\btau)/\kk)
&=\trdeg(\wt M(2\pi i\btau)/\kk)
=\trdeg(\wt M(2\pi i\btau)\point\CC/\CC)
\\
&=\trdeg(F\!\point\CC/\CC)
=\trdeg(\Phi_\kk\point\CC/\CC)
\\
&=\trdeg(\Phi/\CC),
\endalign
$$
which Lemma~1 shows to be equal to $\dim(\SP_{2g}\times Z_g)$.

\proclaim{Step 4}
$M=K^{(2)}$ if $\Gamma=\Gamma_{4,8}$
\ \rom(and consequently\rom, as soon as
$\Gamma\subset\Gamma_{4,8}$\rom).
\endproclaim

We must show that for any modular function~$\lambda_0$ relative
to~$\Gamma_{4,8}$, the components of~$\bdelta^{(3)}\lambda_0$
are rational (rather than just algebraic) over the
field~$K^{(2)}$. In view of Igusa's theorem,
$\lambda_0$~ can be expressed as a rational function
in the quotients $\{\vt_{\ba}/\vt_{\0},\ \ba\in(\ZZ/2\ZZ)^{2g}\}$,
and Proposition~3 below implies that both $\bdelta^{(2)}\lambda_0$
and $\bdelta(\bdelta\vt_{\0}/\vt_{\0})$
are defined over $K^{(1)}(\bdelta\vt_{\0}/\vt_{\0})$.
Consequently,
$\bdelta^{(3)}\lambda_0\in K^{(2)}(\bdelta\vt_{\0}/\vt_{\0})$.
Since $\bdelta\vt_{\0}/\vt_{\0}\in K^{(2)}$ by Step~1,
it follows that $\bdelta^{(3)}\lambda_0\in K^{(2)}$.

\medskip

We conclude this section with two remarks on Theorem~1 and its proof.

\remark{Remark \rom3}
In the case of an arbitrary family of
abelian varieties $f\:\bA\to S$, the algebraic group~$\SP_{2g}$
must be replaced by the {\it Hodge group\/}~$G_A$
of its generic fiber~$A$. Although not always, an equality
$\Galnabla=G_A\otimes\CC$ as in Lemma~1 of~\S\,1 often holds,
cf.~\cite{A1}; the above method can then be extended
to the study of the corresponding field of automorphic functions.
One can thus show that the differential field generated
by Hilbert modular functions (relative to a totally real number field
of degree~$g$) has transcendence degree~$3g$ over~$\CC$.
See~\cite{Re} for some results in this direction,
and for an explicit form of the corresponding differential equations.
\endremark

\remark{Remark \rom4}
Going back to the generic situation where $S=S(\Gamma)$,
denote by~$\Psi$ the field of definition
of the `first periods' of $A/K$, i.e\. the compositum of the field
$$
\Psi_\kk
=\kk\largl(\lambda,\frac1{2\pi i}\Omega_1(\lambda),
\frac1{2\pi i}\Eta_1(\lambda)\largr)
$$
with $\CC$, and view the group~$Z_g$ as an algebraic subgroup
of~$\SP_{2g}$ via the usual map
$\{U\in Z_g\}\mapsto
\bigl\{\bigl(\smallmatrix\1_g&U\\0&\1_g\endsmallmatrix\bigr)
\in\SP_{2g}\bigr\}$.
By Lemma~1 (i.e\. by Picard--Lefschetz theory),
the field~$\Psi$ coincides with the subfield
of the Picard--Vessiot extension $\Phi/\CC(\lambda)$
invariant under~$Z_g$. Combined with the
relation $\xi(\Psi_\kk)=\wt M$ from Remark~2,
this provides a direct proof that
$\trdeg(M/K)=\codim_{\sSP_{2g}}(Z_g)$ when $\kk=\CC$.

When $\kk=\QQ^{\alg}$, this point of view
(extended to the case of Shimura varieties,
as in the previous remark) may prove useful in the study
of the transcendence degree of the field of periods
of abelian varieties defined over~$\QQ^{\alg}$.
In fact, the inclusions
$K^{(1)}\subset\kk(\lambda,P)\subset(K^{(1)})^{\alg}$,
$K^{(2)}\subset\wt M=\xi(\Psi_\kk)\subset(K^{(2)})^{\alg}$
obtained in the course of our proof
(cf.~\S\,3, Remark~2 and \S\,4, Step~2), i.e.
$$
\gathered
\kk(\lambda,\bdelta\lambda)
\subset\kk\largl(\lambda,\frac1{2\pi i}\wt\Omega_1\largr)
\subset\bigl(\kk(\lambda,\bdelta\lambda)\bigr)^{\alg},
\\
\kk(\lambda,\bdelta\lambda,\bdelta^{(2)}\lambda)
\subset\kk\largl(\lambda,\frac1{2\pi i}\wt\Omega_1,
\frac1{2\pi i}\wt\Eta_1\largr)
\subset\bigl(\kk(\lambda,\bdelta\lambda,
\bdelta^{(2)}\lambda)\bigr)^{\alg},
\endgathered
$$
generalize the classical fact that the quotient by~$2\pi i$
of the periods and quasi-periods of an elliptic curve along
a locally invariant cycle can be expressed both as values
of hypergeometric functions and in terms of modular forms.
See~\cite{A2} for an approach to Chudnovsky's theorem on
$\{\omega/2\pi i,\eta/2\pi i\}$ based on these relations,
which, in genus~1, translate into the classical formulae
$$
\gather
\frac{\omega_1}{2\pi}
={}_2\!F_1\largl(\frac12,\frac12,1;\lambda\largr)
=\vt_{00}^2,
\\
\frac{\eta_1}{2\pi}
={}_2\!F_1\largl(-\frac12,\frac12,1;\lambda\largr)
=\frac1{3\vt_{00}^2}\largl(2\vt_{00}^4-\vt_{10}^4
-4\frac1{\pi i}\frac d{d\tau}\Log(\vt_{00}\vt_{10}\vt_{01})\largr),
\endgather
$$
where $\lambda=(\vt_{10}/\vt_{00})^4$~ is
Legendre's modular function, whose derivative satisfies
$$
\frac1{\pi i}\frac{d\lambda}{d\tau}
=\lambda\vt_{01}^4=\lambda(1-\lambda)\vt_{00}^4.
$$
\endremark

\head
\S\,5. Miscellaneous on theta functions
\endhead

In this section, we give the explicit formulae on
derivatives of theta functions already used or mentioned in
Steps~4,~2, and~1 of~\S\,4.

\medskip

We define the (abelian) {\it theta function with characteristic\/}
$\ba=(\ba',\ba'')\in\ZZ^{2g}$ by the convergent series
$$
\vt_{\ba}(\bz)
:=\vt_{\ba}(\bz,\btau)
=\sum_{\bn\in\ZZ^g}
\exp\bigl(\pi i\trans{\bigl(\bn+\half\ba'\bigr)}
\btau\bigl(\bn+\half\ba'\bigr)
+2\pi i\trans{\bigl(\bn+\half\ba'\bigr)}
\bigl(\bz+\half\ba''\bigr)\bigr)
$$
where $\bz\in\CC^g$ is a vector-column and $\btau\in\fH_g$.
The quasi-periodicity of these functions with respect
to the lattice $\ZZ^g+\btau\ZZ^g\subset\CC^g$
allows us to consider only {\it reduced\/} characteristics
$\ba\in\fK=\{0,1\}^{2g}$, $\#\fK=2^{2g}$.
In the customary fashion, we identify $\{0,1\}$ with its
image in~$\ZZ/2\ZZ$, and set for all
$\ba=(\ba',\ba'')$, $\bb=(\bb',\bb'')\in\fK$:
$$
|\ba|=\trans{\ba'}\cdot\ba'',
\qquad
\hbox{$<$}\ba,\bb\hbox{$>$}
=\trans{\ba'}\bb''-\trans{\bb'}\ba''
\equiv|\ba+\bb|+|\ba|+|\bb|\pmod2.
$$

{\it Thetanulls\/} are the values of even theta functions
at the point $\bz=\0$. When no confusion may arise,
we denote them by $\vt_{\ba}=\vt_{\ba}(\btau):=\vt_{\ba}(\0,\btau)$.
Since
$\vt_{\ba}(-\bz,\btau)=(-1)^{|\ba|}\vt_{\ba}(\bz,\btau)$,
a theta function~$\vt_{\ba}(\bz,\btau)$
is even if and only if the number $|\ba|$~ is even.
Hence all non-zero thetanulls are assigned to {\it even\/}
characteristics from the set $\fK_+=\{\ba\in\fK:|\ba|\equiv0\pmod2\}$,
$\#\fK_+=2^{g-1}(2^g+1)$, and trivial ones correspond to {\it odd\/}
characteristics from $\fK_-=\fK\setminus\fK_+$.
We recall (cf.~\cite{I}, pp. 185 and 171) that the thetanulls
are modular forms of weight~$\frac12$ for $\Gamma_{4,8}$.

\medskip

In order to describe the partial logarithmic
derivatives of the thetanulls
$$
\psi_{\ba,jl}=\psi_{\ba,jl}(\btau)
:=\frac{\delta_{jl}\vt_{\ba}}{\vt_{\ba}}=\psi_{\ba,lj},
\qquad \ba\in\fK_+, \quad j,l=1,\dots,g,
\tag1
$$
with respect to the derivations~$\bdelta$:
$$
\delta_{jj}=\frac1{\pi i}\frac\partial{\partial\tau_{jj}},
\quad j=1,\dots,g,
\qquad
\delta_{jl}=\frac1{2\pi i}\frac\partial{\partial\tau_{jl}}
=\delta_{lj},
\quad j,l=1,\dots,g, \; j\ne l,
$$
of~\S\,1, we use the following conventions.
To any meromorphic function $f\:\fH_g\to\CC$,
we assign a meromorphic function with values
in the space of quadratic forms in $\bu\in\CC^g$
by the formula
$$
\bdelta f(\bu)=\sum_{j,l=1}^g\delta_{jl}f\cdot u_ju_l.
$$
Then,
$$
\psi_{\ba}(\bu)=\frac{\bdelta\vt_{\ba}(\bu)}{\vt_{\ba}}
=\sum_{j,l=1}^g\psi_{\ba,jl}\cdot u_ju_l
$$
is the quadratic form corresponding to the symmetric matrix
$\bpsi_{\ba}=(\psi_{\ba,jl})_{j,l=1,\dots,g}$, $\ba\in\fK_+$.
To two quadratic forms $\phi,\eta$ in $\bu\in\CC^g$,
we attach the quartic form $\phi\otimes\eta:=\phi\eta$
given by
$$
\phi\eta(\bu)
=\sum_{j,l,m,p=1}^g\phi_{jl}\eta_{mp}\cdot u_ju_lu_mu_p,
$$
and when $\phi$ has meromorphic coefficients,
we denote by~$\bdelta\phi$ the quartic form
$$
\bdelta\phi(\bu)
=\sum_{j,l,m,p=1}^g\delta_{jl}\phi_{mp}\cdot u_ju_lu_mu_p.
$$

With~\thetag{1} and these conventions in mind,
the system of differential equations on which Step~4 of~\S\,4
is based can then be stated as an equality between
{\it quartic\/} forms, as follows.

\proclaim{Proposition 3 \rm(\cite{Z}, Theorem~1)}
For all $\ba\in\fK_+$\rom,
the thetanulls satisfy the system of second order
partial differential equations
$$
\vt_{\ba}^4\cdot\bdelta\psi_{\ba}
=\frac1{2^{g-2}}\sum_{\bb\in\fK_+}
(-1)^{<\ba,\bb>}\vt_{\bb}^4\cdot\psi_{\bb}^2
-2\vt_{\ba}^4\cdot\psi_{\ba}^2.
$$
\endproclaim

To allow for a comparison with the case of low degrees~$g$
in~\S\,6, consider the rings
$$
Q_g=\QQ[\vt_{\ba},\psi_{\ba,jl}]_{\ba\in\fK_+;j,l=1,\dots,g}
\quad\text{and}\quad
Q_g'=\QQ[\psi_{\ba,jl}]_{\ba\in\fK_+;j,l=1,\dots,g}.
\tag2
$$
Proposition~3 then implies that the {\it fraction field\/} of~$Q_g$
is stable under~$\bdelta$, and (in the notations of~\S\,4) that its
compositum with~$\kk$ coincides with the $\bdelta$-differential field
$M(\vt_{\0})$ generated by $\kk(\vt_{\ba})_{\ba\in\fK_+}=K(\vt_{\0})$.

\medskip

Contrary to the case $g=1$,
where Jacobi's well-known formula expresses the
$z$-derivative at $z=0$ of the unique odd theta function
as a product of thetanulls, the $\bz$-derivatives
of odd theta functions at $\bz=\0$ in higher degrees
are not modular forms. However, they are algebraic over
the field $K^{(1)}$ of~\S\,4, and integral
over the ring~$Q_g$. Indeed:

\proclaim{Proposition 4}
For all $\ba\in\fK_-$\rom, $j=1,\dots,g$\rom,
the following equalities hold\rom:
$$
\align
\text{\rm(i)} &\quad
\largl(\frac1{2\pi i\vt_{\0}}\frac{\partial\vt_{\ba}}{\partial z_j}
\larg|_{\bz=\0}\largr)^2
=\frac1{2^{g-1}}\sum_{\bb\in\fK_+}(-1)^{\trans{\ba'}\bb''}
\largl(\frac{\vt_{\ba+\bb}}{\vt_{\0}}\largr)^2
\largl(\frac{\vt_{\bb}}{\vt_{\0}}\largr)^2
(\psi_{\ba+\bb,jj}-\psi_{\0,jj});
\\
\text{\rm(ii)} &\quad
\largl(\frac1{2\pi i}\frac{\partial\vt_{\ba}}{\partial z_j}
\larg|_{\bz=\0}\largr)^4
=\frac1{2^{g-2}}\sum_{\bb\in\fK_+}
(-1)^{|\ba+\bb|}\vt_{\bb}^4\psi_{\bb,jj}^2.
\endalign
$$
\endproclaim

\demo{Proof}
(i) We shall use the quartic Riemann relations
in the following form, valid for {\it all\/} characteristics
$\ba,\bc\in\fK$:
$$
\vt_{\ba+\bc}^2(\bz)\vt_{\ba}^2(\bz)
=\frac1{2^g}\sum_{\bb\in\fK}
(-1)^{<\ba,\bb>}(-1)^{\trans{\bc'}(\ba''+\bb'')}
\vt_{\bb+\bc}(2\bz)\vt_{\bb+\bc}(\0)\vt_{\bb}^2(\0)
\tag3
$$
(see~\cite{Kr}, VII, \S\,10),
together with the heat equation
$$
\delta_{jl}\vt_{\ba}(\bz,\btau)
=\frac1{(2\pi i)^2}
\frac{\partial^2}{\partial z_j\,\partial z_l}\vt_{\ba}(\bz,\btau),
\qquad \ba\in\fK, \quad j,l=1,\dots,g
$$
(see, e.g., \cite{Kr}, I, \S\,5, Satz~13), which
allows to write the Taylor expansions
of {\it even\/} theta functions as
$$
\vt_{\bb}(\bz)
=\vt_{\bb}\cdot\largl(1+\frac{(2\pi i)^2}2\psi_{\bb}(\bz)\largr)
+O(\bz^4),
\qquad \bb\in\fK_+.
$$
Notice that in Riemann's formula~\thetag{3}, only even
characteristics
$\bb,\bb+\bc$ actually appear on the right-hand side.
Setting
$\ba=\bc\in\fK_-$ in~\thetag{3}, we derive
$$
\vt_{\0}^2(\bz)\vt_{\ba}^2(\bz)
=\frac1{2^g}\sum_{\bb\in\fK_+}(-1)^{\trans{\ba'}\bb''}
\vt_{\ba+\bb}(2\bz)\vt_{\ba+\bb}(\0)\vt_{\bb}^2(\0).
$$
Since $\ba$ is odd, the Taylor expansion of the left-hand
side of this formula reads
$$
\vt_{\0}^2\cdot\largl(\trans{\bz}
\cdot\frac{\partial\vt_{\ba}}{\partial\bz}\larg|_{\bz=\0}\largr)^2
+O(\bz^4).
$$
Equating quadratic terms, we obtain
$$
\vt_{\0}^2\cdot\largl(\trans{\bz}
\cdot\frac{\partial\vt_{\ba}}{\partial\bz}\larg|_{\bz=\0}\largr)^2
=\frac{(2\pi i)^2}{2^{g-1}}\sum_{\bb\in\fK_+}
(-1)^{\trans{\ba'}\bb''}\vt_{\ba+\bb}^2\vt_{\bb}^2
\cdot\psi_{\ba+\bb}(\bz).
$$
Furthermore,
$$
\sum_{\bb\in\fK_+}(-1)^{\trans{\ba'}\bb''}
\vt_{\ba+\bb}^2\vt_{\bb}^2=0,
$$
since no constant term appears on the left-hand side.
Multiplying the latter relation by $\psi_{\0}(\bz)$
and substracting, we finally get
$$
\vt_{\0}^2\cdot\largl(\trans{\bz}
\cdot\frac{\partial\vt_{\ba}}{\partial\bz}\larg|_{\bz=\0}\largr)^2
=\frac{(2\pi i)^2}{2^{g-1}}\sum_{\bb\in\fK_+}
(-1)^{\trans{\ba'}\bb''}\vt_{\ba+\bb}^2\vt_{\bb}^2
\cdot(\psi_{\ba+\bb}-\psi_{\0})(\bz),
$$
and Proposition~4~(i) follows on considering
the $z_j^2$-term of these quadratic forms.

(ii) Setting $\ba\in\fK_-$, $\bc=\0$
in Riemann's formula~\thetag{3},
and developing the Taylor expansions
$$
\vt_{\bb}(\bz)
=\vt_{\bb}\cdot\largl(1+\frac{(2\pi i)^2}2\psi_{\bb}(\bz)
+\frac{(2\pi i)^4}{24}(\psi_{\bb}^2+\bdelta\psi_{\bb})(\bz)\largr)
+O(\bz^6),
\qquad \bb\in\fK_+,
$$
of even theta functions to the fourth order
(cf.~\cite{Z}, \S\,2, Lemma~2), we derive
from Formula~\thetag{3} the equality of quartic forms
$$
\largl(\trans{\bz}
\cdot\frac{\partial\vt_{\ba}}{\partial\bz}\larg|_{\bz=\0}\largr)^4
=\frac{(2\pi i)^4}{3\cdot2^{g-1}}\sum_{\bb\in\fK_+}
(-1)^{<\ba,\bb>}\vt_{\bb}^4\cdot(\psi_{\bb}^2+\bdelta\psi_{\bb})(\bz).
$$
In view of Proposition~3 and the identity
$$
\sum_{\bb\in\fK_+}(-1)^{<\ba+\bc,\bb>}
=(-1)^{|\ba+\bc|}\cdot2^{g-1},
\qquad \ba\in\fK_-, \quad \bc\in\fK_+
$$
(cf.~\cite{Z}, \S\,3, Lemma~7), this transforms into
$$
\largl(\trans{\bz}
\cdot\frac{\partial\vt_{\ba}}{\partial\bz}\larg|_{\bz=\0}\largr)^4
=\frac{(2\pi i)^4}{2^{g-2}}\sum_{\bb\in\fK_+}
(-1)^{|\ba+\bb|}\vt_{\bb}^4\cdot\psi_{\bb}^2(\bz),
\qquad \ba\in\fK_-,
$$
and Proposition 4~(ii) follows on considering the
$z_j^4$-term of these quartic forms.
\enddemo

The next statement holds for all pairs
$\ba,\bb\in\fK_+$ of distinct even characteristics (and a
sharper result will be given in Formula~\thetag{6a} below in
the case of degree~$2$), but one instance suffices for our
purpose (cf.~\S\,4, Step~1). Recall that $\bpsi_{\ba}$
denotes the symmetric matrix $\bdelta\vt_{\ba}/\vt_{\ba}$
attached to the quadratic form~$\psi_\ba$.

\proclaim{Lemma 4}
There exists even characteristics $\ba,\bb$ such that
$\det(\bpsi_{\bb}-\bpsi_{\ba})$ is a non-zero meromorphic
modular form of weight~$2$ with respect to~$\Gamma_{4,8}$.
\endproclaim

\demo{Proof}
Since $\lambda_0=\vt_{\bb}/\vt_{\ba}$ is a modular function,
$$
(\bdelta\lambda_0)(\gamma\btau)
=(c\btau+d)\cdot(\bdelta\lambda_0)(\btau)\cdot\trans{(c\btau+d)},
\qquad
\gamma=\bigl(\smallmatrix a&b\\c&d\endsmallmatrix\bigr)
\in\Gamma_{4,8}
\tag4
$$
(see~\cite{Z}, \S\,8, Lemma~16),
where $\bdelta\lambda_0$ is arranged as a symmetric matrix,
and we obtain that
$$
\eta_{\ba,\bb}
:=\det\largl(\frac{\bdelta\lambda_0}{\lambda_0}\largr)
=\det(\bpsi_{\bb}-\bpsi_{\ba}),
\qquad \ba,\bb\in\fK_+, \quad \ba\ne\bb,
$$
is a meromorphic modular form of weight~$2$.

We now show that $\eta_{\ba,\bb}$ is non-zero for
$\ba=(\ba',\ba'')=(\0,\0)\in\fK_+$ and
$\bb=\allowmathbreak(\bb',\bb'')=(\1,\0)\in\fK_+$,
where $\1$~ is a $g$-vector-column with unit entries.
Choosing $\btau=\tau\1_g\in\fH_g$, where $\tau\in\fH_1$,
we get from the definition of the theta functions
$$
\vt_{\ba}(\bz,\tau\1_g)
=\prod_{j=1}^g\vt_{00}(z_j,\tau)
\quad\text{and}\quad
\vt_{\bb}(\bz,\tau\1_g)
=\prod_{j=1}^g\vt_{10}(z_j,\tau);
$$
here, $\vt_{00}(z,\tau)$ and $\vt_{10}(z,\tau)$
stand for the usual elliptic theta functions.
Taking $\partial^2/\partial z_j\,\partial z_l$-derivatives
and evaluating at $\bz=\0$, we deduce from the heat equations
in degrees~$g$ and~$1$, and from the even
character of~$\vt_{00}(z,\tau)$, that
$$
\psi_{\ba,jl}(\tau\1_g)=\cases
\dsize\frac1{\vt_{00}}
\frac{\partial\vt_{00}}{\pi i\,\partial\tau}(0,\tau)
:=\psi_{00}(\tau) & \text{if $j=l$}, \\
0 & \text{if $j\ne l$}.
\endcases
$$
Similarly,
$\psi_{\bb,jl}(\tau\1_g)$ vanishes if $j\ne l$
and is otherwise equal to $\psi_{10}(\tau)$.
Thus, $(\bpsi_{\bb}-\bpsi_{\ba})(\tau\1_g)$
is a diagonal matrix, with determinant
$$
\eta_{\ba,\bb}(\tau\1_g)
=\det(\bpsi_{\bb}-\bpsi_{\ba})(\tau\1_g)
=\bigl(\psi_{10}(\tau)-\psi_{00}(\tau)\bigr)^g.
$$
Now, Legendre's modular function $\lambda=(\vt_{10}/\vt_{00})^4$
is not constant, so that
$$
\psi_{10}-\psi_{00}
=\frac1{4\lambda}\frac{\partial\lambda}{\pi i\,\partial\tau}
\ne0,
$$
and $\eta_{\ba,\bb}$ is indeed a non-zero degree~$g$
modular form.
\enddemo

\head
\S\,6. The differential ring of thetanulls in genus~$2$
\endhead

Theorem~2 can be sharpened when $g=1$ and $g=2$,
because in both of these cases, the {\it rings\/}
$$
Q_g=\QQ[\vt_{\ba},\psi_{\ba,jl}]_{\ba\in\fK_+;j,l=1,\dots,g}
\quad\text{and}\quad
Q_g'=\QQ[\psi_{\ba,jl}]_{\ba\in\fK_+;j,l=1,\dots,g}.
$$
generated over~$\QQ$ by the thetanulls and their
logarithmic derivatives (cf.~\S\,5, \thetag{1},~\thetag{2})
are themselves stable under derivation.

\medskip

First, we consider the well-known case $g=1$
where we have three even characteristics,
$\fK_+=\{00,01,10\}$. In this case there exists only one
parameter $\tau=\tau_{11}$, and the notations
$\delta$, $\psi_{\ba}$, need no $jl$-indexation. The
rings $Q_1$ and $Q_1'$ are $\delta$-stable, since
$\delta\vt_{\ba}=\vt_{\ba}\psi_{\ba}$, $\ba\in\fK_+$, by definition,
while
$$
\aligned
\delta\psi_{10}&=2(\psi_{10}\psi_{00}+\psi_{10}\psi_{01}-\psi_{00}\psi_{01}),
\\
\delta\psi_{00}&=2(\psi_{10}\psi_{00}+\psi_{00}\psi_{01}-\psi_{10}\psi_{01}),
\\
\delta\psi_{01}&=2(\psi_{10}\psi_{01}+\psi_{00}\psi_{01}-\psi_{10}\psi_{00})
\endaligned
$$
(this system was discovered by G.~Halphen in~1881).
Moreover, $Q_1$~ is integral over~$Q_1'$, in view of the formulae
$$
\vt_{00}^4=4(\psi_{10}-\psi_{01}),
\qquad
\vt_{01}^4=4(\psi_{10}-\psi_{00}),
\qquad
\vt_{10}^4=4(\psi_{00}-\psi_{01})
$$
(see \cite{Z}, Introduction).
By~\cite{Ra},~\cite{M} or Theorem~2, both rings $Q_1$
and~$Q_1'$ have transcendence degree~$3$ over~$\QQ$.

\medskip

In the case $g=2$, we use for simplicity the map
$(\ZZ/2\ZZ)^2\to\{0,1,2,3\}$,
$$
(0,0)\mapsto0, \quad (0,1)\mapsto1, \quad
(1,0)\mapsto2, \quad (1,1)\mapsto3,
$$
to represent a characteristic
$\ba=(\ba',\ba'')\in(\ZZ/2\ZZ)^2\times(\ZZ/2\ZZ)^2$
by two digits only. Then,
$$
\fK_+=\{00,01,02,03,10,12,20,21,30,33\}.
$$
We renumerate the entries of the matrix~$\btau$ as
$\tau_1=\tau_{11}$, $\tau_2=\tau_{22}$, and $\tau_3=\tau_{12}$,
and proceed similarly with the derivations
$$
\bdelta=\largl\{
\delta_1=\frac1{\pi i}\frac\partial{\partial\tau_1}, \
\delta_2=\frac1{\pi i}\frac\partial{\partial\tau_2}, \
\delta_3=\frac1{2\pi i}\frac\partial{\partial\tau_3}
\largr\},
$$
and with the $\psi_{\ba}$-notation.

\proclaim{Theorem 3}
\rom{(i)} In the case $g=2$\rom, the rings
$$
Q_2=\QQ[\vt_{\ba},\psi_{\ba,j}]_{\ba\in\fK_+;j=1,2,3}
\quad\text{and}\quad
Q_2'=\QQ[\psi_{\ba,j}]_{\ba\in\fK_+;j=1,2,3}
$$
are stable under the derivations~$\delta_j$\rom, $j=1,2,3$.

\rom{(ii)} All thetanulls are algebraic over
$\QQ(\psi_{\ba,j})_{\ba\in\fK_+;j=1,2,3}$.

\rom{(iii)} Both rings $Q_2$ and~$Q_2'$ have transcendence
degree~$10$ over~$\QQ$.

\rom{(iv)} A possible choice for ten elements of~$Q_2$\rom,
algebraically independent over~$\QQ$\rom, is given by
$$
\vt_{00},\vt_{01},\vt_{02}, \
\psi_{00,1},\psi_{01,1},\psi_{02,1}, \
\psi_{00,2},\psi_{01,2},\psi_{02,2}, \
\psi_{00,3}.
$$
\endproclaim

\demo{Proof}
(i) Let $\ba_1=\ba$, $\ba_2=\ba+\bc$,
$\ba_3=\ba+\bd$, $\ba_4=\ba+\bc+\bd$ be four different even
characteristics (such a collection $\{\ba_1,\ba_2,\ba_3,\ba_4\}$
is called a {\it G\"opel system}, and there exist
fifteen G\"opel systems). Then
$$
\bdelta(\psi_{\ba_1}+\psi_{\ba_2}+\psi_{\ba_3}+\psi_{\ba_4})
=(\psi_{\ba_1}+\psi_{\ba_2}+\psi_{\ba_3}+\psi_{\ba_4})^2
-2(\psi_{\ba_1}^2+\psi_{\ba_2}^2+\psi_{\ba_3}^2+\psi_{\ba_4}^2).
\tag5a
$$
System~\thetag{5a}, which was obtained in~\cite{O},
provides an expression for each $\delta\psi_{\ba}$, $\ba\in\fK_+$.
Namely (cf.~\cite{Z}, \S\,6, Formulae~(6.17)):
$$
\bdelta\psi_{\ba}
=-2\psi_{\ba}^2
-\frac13\sum_{\bb\in\fK_+}\psi_{\bb}^2
-\frac16\largl(\sum_{\bb\in\fK_+}\psi_{\bb}\largr)^2
+\frac14\sum_{G\owns\ba}\largl(\sum_{\bb\in G}\psi_{\bb}\largr)^2,
\qquad \ba\in\fK_+,
\tag5b
$$
where the summation $\sum_{G\owns\ba}$ goes over all different
G\"opel systems containing~$\ba$. The differential stability
of $Q_2'$ is clear from~\thetag{5b}, and that of $Q_2$ then
follows from the defining equations~\thetag{1}.

(ii) This result (which is not the statement of
algebraicity over~$K^{(1)}$ given in Step~1 of~\S\,4) is
shown in Theorem~6 of~\cite{Z}. Here is a sketch of its
proof. Let $\ba,\bb\in\fK_+$, $\ba\ne\bb$.
As noticed in~\S\,5, Lemma~4,
$$
\eta_{\ba,\bb}=\eta_{\bb,\ba}
:=\det(\bpsi_{\ba}-\bpsi_{\bb})
=(\psi_{\ba,1}-\psi_{\bb,1})(\psi_{\ba,2}-\psi_{\bb,2})
-(\psi_{\ba,3}-\psi_{\bb,3})^2
$$
is a meromorphic modular form of weight~$2$, so that
$\vt_{\ba}^2\vt_{\bb}^2\eta_{\ba,\bb}$ is a
holomorphic modular form of weight~$4$
with respect to~$\Gamma_{4,8}$. By the sharper version
of Igusa's theorem in the genus~$2$ case, any such
form can be expressed as a polynomial in thetanulls, and indeed,
$$
\eta_{\ba,\bb}
=\pm\frac1{2^4}\prod_{\bc\in\fK_+}\vt_{\bc}^2
\cdot\prod_{G\supset\{\ba,\bb\}}
\largl(\prod_{\bd\in G}\vt_{\bd}^{-2}\largr),
\tag6a
$$
(see \cite{Z}, \S\,8, Lemma~20~(a), and
Formula~\thetag{6c} below for an explicit expression);
here, the product $\prod_{G\supset\{\ba,\bb\}}$ goes over all
(i.e\. the two) G\"opel systems containing $\ba$ and~$\bb$,
and the only term that appears in the denominator of the
right-hand side of~\thetag{6a} is $\vt_{\ba}^2\vt_{\bb}^2$.

Multiplying Formulae~\thetag{6a} for fixed~$\ba$ over all $\bb\ne\ba$,
and then over all different pairs $\ba,\bb$, we get
(cf.~\cite{Z},~(8.30))
$$
\vt_{\ba}^{72}
=\pm2^{-4\cdot27}\prod_{\bc\in\fK_+}\vt_{\bc}^{18}
\cdot\prod\Sb\bb\in\fK_+\\\bb\ne\ba\endSb\eta_{\ba,\bb}^{-3}
=\pm2^{72}
\prod\Sb\text{different}\\\text{pairs $\bc,\bd$}\endSb\eta_{\bc,\bd}
\cdot\prod\Sb\bb\in\fK_+\\\bb\ne\ba\endSb\eta_{\ba,\bb}^{-3},
\qquad \ba\in\fK_+.
\tag6b
$$
Consequently, all thetanulls are algebraic
over $\QQ(\psi_{\ba,j})_{\ba\in\fK_+;j=1,2,3}$.

(iii) Recalling the notations of~\S\,4, we infer from~(i)
and (ii) above that $M$ and the fraction fields of
$Q_2'$ and of $Q_2$ have the same algebraic closure. Since
$2g^2+g=10$, the result immediately follows from our Theorem~2.

(iv) We shall find a maximal set of algebraically independent elements
in the ring~$Q_2$ by exhibiting thirty independent relations
between its fourty generators.%
\footnote[5]{We do not claim that these relations
generate a full ideal of definition for~$Q_2$.}

Specializing the Riemann relations~\thetag{3} used in the
proof of Proposition~4 to the case
$\bc\allowmathbreak\in\{01,02,03\}$, $\ba=00$, we obtain
$$
\vt_{00}^2\vt_{01}^2-\vt_{02}^2\vt_{03}^2
=\vt_{20}^2\vt_{21}^2,
\quad
\vt_{00}^2\vt_{02}^2-\vt_{01}^2\vt_{03}^2
\vt_{10}^2\vt_{12}^2,
\quad
\vt_{00}^2\vt_{03}^2-\vt_{01}^2\vt_{02}^2
=\vt_{30}^2\vt_{33}^2,
\tag3a
$$
while the case $\bc=00$ yields
$$
\vt_{00}^4-\vt_{01}^4=\vt_{10}^4+\vt_{33}^4,
\quad
\vt_{00}^4-\vt_{02}^4=\vt_{21}^4+\vt_{30}^4,
\quad
\vt_{00}^4-\vt_{03}^4=\vt_{12}^4+\vt_{20}^4.
\tag3b
$$
In view of~\thetag{3a} and \thetag{3b},
all thetanulls can be expressed in terms of
the four functions
$\{\vt_{00},\vt_{01},\vt_{02},\allowmathbreak\vt_{03}\}$,
which must therefore be algebraically
independent. We shall denote by $\fK_0:=\{00,01,02,03\}$
the corresponding set of characteristics.

\smallskip

By $\bdelta$-differention of~\thetag{3a} and~\thetag{3b},
we obtain each $\psi_{\ba,j}$, $\ba\in\fK_+$, $j=1,2,3$,
as a linear combination with coefficients in
$\QQ(\vt_{\ba})_{\ba\in\fK_+}$
of the functions $\psi_{\ba,j}$, $\ba\in\fK_0$, $j=1,2,3$.
We are thus reduced to finding six independent relations linking
the thirteen functions $\vt_{03}$ and $\psi_{\ba,j}$,
$\ba\in\fK_0$, $j=1,2,3$, over
$\QQ(\vt_{00},\vt_{01},\vt_{02})$.

\smallskip

Developing the functions
$\eta_{\ba,\bb}$, $\ba\ne\bb\in\fK_0$,
used in~(ii) above and expliciting the right-hand side
of Formula~\thetag{6a}, we obtain:
$$
\aligned
\eta_{00,01}
:=(\psi_{00,1}-\psi_{01,1})(\psi_{00,2}-\psi_{01,2})
-(\psi_{00,3}-\psi_{01,3})^2
&=\phantom+\frac1{16}
\frac{\vt_{10}^2\vt_{12}^2\vt_{30}^2\vt_{33}^2}{\vt_{00}^2\vt_{01}^2},
\\
\eta_{00,02}
:=(\psi_{00,1}-\psi_{02,1})(\psi_{00,2}-\psi_{02,2})
-(\psi_{00,3}-\psi_{02,3})^2
&=\phantom+\frac1{16}
\frac{\vt_{20}^2\vt_{21}^2\vt_{30}^2\vt_{33}^2}{\vt_{00}^2\vt_{02}^2},
\\
\eta_{01,02}
:= (\psi_{01,1}-\psi_{02,1})(\psi_{01,2}-\psi_{02,2})
-(\psi_{01,3}-\psi_{02,3})^2
&=-\frac1{16}
\frac{\vt_{10}^2\vt_{12}^2\vt_{20}^2\vt_{21}^2}{\vt_{01}^2\vt_{02}^2},
\\
\eta_{00,03}
:=(\psi_{00,1}-\psi_{03,1})(\psi_{00,2}-\psi_{03,2})
-(\psi_{00,3}-\psi_{03,3})^2
&=\phantom+\frac1{16}
\frac{\vt_{10}^2\vt_{12}^2\vt_{20}^2\vt_{21}^2}{\vt_{00}^2\vt_{03}^2},
\\
\eta_{01,03}
:=(\psi_{01,1}-\psi_{03,1})(\psi_{01,2}-\psi_{03,2})
-(\psi_{01,3}-\psi_{03,3})^2
&=-\frac1{16}
\frac{\vt_{20}^2\vt_{21}^2\vt_{30}^2\vt_{33}^2}{\vt_{01}^2\vt_{03}^2},
\\
\eta_{02,03}
:=(\psi_{02,1}-\psi_{03,1})(\psi_{02,2}-\psi_{03,2})
-(\psi_{02,3}-\psi_{03,3})^2
&=-\frac1{16}
\frac{\vt_{10}^2\vt_{12}^2\vt_{30}^2\vt_{33}^2}{\vt_{02}^2\vt_{03}^2}.
\endaligned
\tag6c
$$
Notice that in view of Relations~\thetag{3a}, the numerators
of the right-hand terms of~\thetag{6c} are polynomials
in~$\vt_{\ba}$, $\ba\in\fK_0$. For instance, the first one reads
$$
\vt_{10}^2\vt_{12}^2\vt_{30}^2\vt_{33}^2
=(\vt_{00}^2\vt_{02}^2-\vt_{01}^2\vt_{03}^2)
(\vt_{00}^2\vt_{03}^2-\vt_{01}^2\vt_{02}^2)
=-\vt_{00}^2\vt_{01}^2\vt_{03}^4+\dotsb.
$$
Now, consider the expressions
$$
\Chi_1=(\psi_{00,3}-\psi_{01,3})^2, \quad
\Chi_2=(\psi_{00,3}-\psi_{02,3})^2, \quad
\Chi_3=(\psi_{02,3}-\psi_{01,3})^2.
$$
They formally satisfy
$$
\chi_1^2+\chi_2^2+\chi_3^2
-2\Chi_1\Chi_2-2\Chi_2\Chi_3-2\Chi_3\Chi_1=0.
$$
Substituting into this relation the expressions for the
$\chi$'s given by the first three lines of~\thetag{6c} and
developing the numerators of the right-hand terms
of~\thetag{6c} as above, we obtain a polynomial relation
$\fR_0$ between
$$
\vt_{03}, \;
\psi_{00,1},\psi_{00,2}, \;
\psi_{01,1},\psi_{01,2}, \;
\psi_{02,1},\psi_{02,2}
$$
with coefficients in $\QQ(\vt_{00},\vt_{01},\vt_{02})$.
The coefficient of $\vt_{03}^8$ in~$\fR_0$, viewed as a polynomial in
$\vt_{03}$, is the non-zero constant
$-3\cdot\frac1{16^2}$, and we deduce that
$\vt_{03}$~ {\it is algebraic over the field}
$$
\QQ(\vt_{\ba},\psi_{\ba,j};\ \ba\in\fK',\ j=1,2)
\qquad \text{where} \quad
\fK':=\{00,01,02\}=\fK_0\setminus\{03\}.
\tag7
$$
By~\thetag{6c}, the four functions
$\psi_{\ba,3}$, $\ba\in\fK_0$, are then algebraic over
the field generated by anyone of them, say~$\psi_{00,3}$,
over the field~\thetag{7}.

\smallskip

We conclude the proof of Theorem~3~(iv) by finding
(in a similar, though more subtle, way) two
independent relations $\fR_1,\fR_2$ linking
$\psi_{03,1},\psi_{03,2}$ to the nine generators
of the field~\thetag{7}. Consider the expressions
$$
\aligned
\phi_0
&=(\psi_{03,1}-\psi_{01,1})(\psi_{01,1}-\psi_{02,1})\eta_{02,03}
+(\psi_{02,1}-\psi_{03,1})(\psi_{01,1}-\psi_{02,1})\eta_{03,01}
\\ &\qquad
+(\psi_{02,1}-\psi_{03,1})(\psi_{03,1}-\psi_{01,1})\eta_{01,02},
\\
\phi_1
&=(\psi_{03,1}-\psi_{00,1})(\psi_{00,1}-\psi_{02,1})\eta_{02,03}
+(\psi_{02,1}-\psi_{03,1})(\psi_{00,1}-\psi_{02,1})\eta_{03,00}
\\ &\qquad
+(\psi_{02,1}-\psi_{03,1})(\psi_{03,1}-\psi_{00,1})\eta_{00,02},
\\
\phi_2
&=(\psi_{03,1}-\psi_{00,1})(\psi_{00,1}-\psi_{01,1})\eta_{01,03}
+(\psi_{01,1}-\psi_{03,1})(\psi_{00,1}-\psi_{01,1})\eta_{03,00}
\\ &\qquad
+(\psi_{01,1}-\psi_{03,1})(\psi_{03,1}-\psi_{00,1})\eta_{00,01},
\\
\phi_3
&=(\psi_{02,1}-\psi_{00,1})(\psi_{00,1}-\psi_{01,1})\eta_{01,02}
+(\psi_{01,1}-\psi_{02,1})(\psi_{00,1}-\psi_{01,1})\eta_{02,00}
\\ &\qquad
+(\psi_{01,1}-\psi_{02,1})(\psi_{02,1}-\psi_{00,1})\eta_{00,01}.
\endaligned
\tag8a
$$
When the $\eta$'s are developed as in the middle terms
of~\thetag{6c}, a computation shows that they formally
satisfy the relation
$$
\bigl((\phi_0+\phi_1+\phi_2+\phi_3)^2
-2(\phi_0^2+\phi_1^2+\phi_2^2+\phi_3^2)\bigr)^2
-64\phi_0\phi_1\phi_2\phi_3=0.
\tag8b
$$
Substituting into~\thetag{8a} the expressions for the
$\eta$'s given by
the right-hand terms of~\thetag{6c}, we deduce
from~\thetag{8b} a polynomial relation $\fR_1$ between
$$
\text{$\psi_{00,1}$, $\psi_{01,1}$,
$\psi_{02,1}$, $\psi_{03,1}$}
$$
with coefficients in $\QQ(\vt_{\ba})_{\ba\in\fK_0}$.
Now, the coefficient of $\psi_{03,1}^8$ in $\fR_1$, viewed as a polynomial in
$\psi_{03,1}$, is
$$
\bigl((\eta_{00,01}+\eta_{00,02}+\eta_{01,02})^2
-2(\eta_{00,01}^2+\eta_{00,02}^2+\eta_{01,02}^2)\bigr)^2,
$$
and we can check that this expression is not identically~$0$
by evaluating it at $\btau=\tau\1_2$ with
$\tau\in\fH_1$, as in the proof of~Lemma~4:
one finds that it reduces to the non-zero genus 1 modular form
$\eta_{01,02}^4(\tau\1_2)
=\frac1{16^4}\vt_{10}^{32}(\tau).
$
We thus deduce from~$\fR_1$ that $\psi_{03,1}$~ is
algebraic over the field generated by
$\vt_{\ba}$, $\psi_{\ba,1}$, $\ba\in\fK'$, and~$\vt_{03}$,
hence over the field~\thetag{7}, in view of the algebraicity
of~$\vt_{03}$. By the same argument,
$\psi_{03,2}$~ is algebraic over the field generated by
$\vt_{\ba}$, $\psi_{\ba,2}$, $\ba\in\fK'$, and $\vt_{03}$,
hence over the field~\thetag{7}. Therefore, the ten functions
$$
\vt_{00},\vt_{01},\vt_{02}, \
\psi_{00,1},\psi_{01,1},\psi_{02,1}, \
\psi_{00,2},\psi_{01,2},\psi_{02,2}, \
\psi_{00,3}
$$
are algebraically independent over~$\QQ$,
and the proof of Theorem~3 is completed.
\enddemo

\goodbreak

\remark{Remark \rom5}
In the same vein, it can be proved in the
case $g=3$ that the {\it fraction field\/} of the ring
$$
Q_3'=\QQ[\psi_{\ba,jl}]_{\ba\in\fK_+;j,l=1,2,3}
$$
is $\bdelta$-stable, and that all thetanulls
are algebraic over it. This follows from~\cite{Z},
Theorems~4 and~6. Theorem~2 of the present paper then
implies that $Q_3'$ has transcendence degree~$21$ over~$\QQ$.
\endremark

\subhead
Acknowledgements
\endsubhead
This work was made under an INTAS--RFBR program,
Contract 97-1904. The second author expresses his
gratitude to the Ostrowski Foundation for financial support.


\enddocument
\end